\documentstyle[amssymb,12pt]{article}

\input{tcilatex}

\begin{document}

\title{Formal Frobenius structures generated by geometric deformation algebras}
\author{Mircea Cr\^a\c sm\u areanu}
\date{}
\maketitle

\begin{abstract}
Necessary and sufficient conditions for some deformation algebras to provide
formal Frobenius structures are given. Also, examples of formal Frobenius
structures with fundamental tensor that is not of the deformation type and
examples of symmetric non-metric connections are presented.
\end{abstract}



\vspace{.3cm}

\section*{Introduction}

Frobenius manifolds were introduced by Dubrovin(\cite{d:b}) as a
''coordinate-free'' framework for Gromov-Witten invariants. Also, Frobenius
manifolds provide a natural geometric setting for understand the
bi-Hamiltonian structure of hydrodynamics systems. The present paper is
devoted to a notion related to Frobenius manifolds, namely {\it formal} {\it %
Frobenius structure}.

Our aim is to connect the geometry of formal Frobenius structures to
classical differential geometry through deformation algebras of connections,
a notion introduced by Izu Vaisman in \cite{i:v}. The starting point of our
study is the following remark: a formal Frobenius structure is a pair $%
\left( g,A\right) $ with $g$ a Riemann metric and $A$ a tensor field of $%
\left( 1,2\right) $-type, subject to the conditions below. Naturally
associated to $g$ is the Levi-Civita connection $\nabla $ and then the pair $%
\left( \nabla ,A\right) $ yields another linear connection $\overline{\nabla
}=\nabla +A$. We search the converse, namely starting with $\left( \nabla ,%
\overline{\nabla }\right) $ we find conditions for $A=\overline{\nabla }%
-\nabla $ to satisfy the definition of a Frobenius structure. Therefore, in
a general sense, every Frobenius structure can be viewed as a geometric
deformation algebra!

An usual condition in the theory of Frobenius manifolds is the symmetry of $%
A $ which is equivalent in our framework with the equality of torsions for $%
\nabla $ and $\overline{\nabla }$. Because $\nabla $ is without torsion we
arrive at condition of torsionless of $\overline{\nabla }$ which restricts
the number of remarkable deformation algebras. So, we add a version of
Frobenius structures without commutativity condition.

Let us sketch the contents: after a first section which reviews the main
definitions, in the following two sections several examples of formal and
weak Frobenius structures provided by deformation algebras are discussed.
Another section is devoted to some examples of Frobenius structures with $A$
not of $\overline{\nabla }-\nabla $ type. Because in the above discussion
the search of torsion-free linear connections appears as very important, the
last section gives examples of symmetric connections which are not
Levi-Civita connection for any Riemannian metric!

{\bf Acknowledgments} I would like to thank Larry Bates and Alexey
Tsygvintsev for sending their papers and Vladimir B\u alan for
useful remarks.

\section{Frobenius structures and deformation algebras}

Let $\left( M,g\right) $ be a fixed Riemannian manifold for which we denote $%
C^{\infty }\left( M\right) $ the ring of smooth real functions, ${\cal X}%
\left( M\right) $ the Lie algebra of vector fields, $T_{s}^{r}\left(
M\right) $ the $C^{\infty }\left( M\right) $-module of tensor fields of $%
\left( r,s\right) $-type, $\Omega ^{k}\left( M\right) $ the $C^{\infty
}\left( M\right) $-module of $k$-differential forms. Let $n=\dim M$ be
finite.

{\bf Definition 1.1}(\cite[p. 286]{s:f}) The triple $\left( M,g,A\right) $
with $A\in T_{2}^{1}\left( M\right) $ is called a {\it formal Frobenius
structure} if: \newline
(i) for every $X,Y,Z\in {\cal X}\left( M\right) $:
$$
g\left( A\left( X,Y\right) ,Z\right) =g\left( X,A\left( Y,Z\right) \right) %
\eqno\left( 1.1\right)
$$
(ii) $A$ is commutative i.e. $A\left( X,Y\right) =A\left( Y,X\right) $.

Using the symmetry of $g$ it follows that $\left( 1.1\right) $ means the
invariance of application $g\left( A\left( \cdot ,\cdot \right) ,\cdot
\right) :{\cal X}\left( M\right) ^{3}\rightarrow {\cal X}\left( M\right) $
with respect to cyclic permutations.

Recall that for the given $A\in T_{2}^{1}\left( M\right) $ the
multiplication $X*Y:=A\left( X,Y\right) $ defines a $C^{\infty }\left(
M\right) $-algebra structure on ${\cal X}\left( M\right) $. Sometimes in the
definition of formal Frobenius structures one also asks for a third
condition, namely the existence of a unit element in this algebra but we do
not work in this framework. For other several types of Frobenius structures
see \cite{h:c}, \cite{m:y}, \cite{s:f2}.

Given two linear connections $\nabla ,\overline{\nabla }$ the $C^{\infty
}\left( M\right) $-algebra defined by $A=\overline{\nabla }-\nabla $ is
called {\it the deformation algebra} in \cite[p. 83]{i:v}( see also \cite
{b:v}). In the cited paper it is proved that the deformation algebra is
commutative if and only if $\nabla $ and $\overline{\nabla }$ have the same
torsion. If we start with $\nabla $ the Levi-Civita connection of $g$ it
follows that $\overline{\nabla }$ must be symmetric (torsionless). This
condition is very restrictive and then, in order to use some remarkable
deformation algebras which are not commutative, we consider:

{\bf Definition 1.2} A triple satisfying only condition (i) of Definition
1.1 is called a {\it weak Frobenius structure}.

\section{Frobenius structures generated by deformation algebras}

{\bf Example 2.1 Subgeodesic correspondences}

Two Riemannian metrics $g,\overline{g}$ are said {\it in a }$g${\it %
-subgeodesic correspondence} if there exists $\theta \in \Omega ^{1}\left(
M\right) $ and $P\in {\cal X}\left( M\right) $ such that $A=\overline{\nabla
}-\nabla $ is given by:
$$
A=\theta \otimes \delta +\delta \otimes \theta +g\otimes P\eqno\left(
2.1\right)
$$
where $\nabla ,\overline{\nabla }$ is the Levi-Civita connection of $g,%
\overline{g}$ and $\delta $ is the Kronecker tensor. Let $\psi \in \Omega
^{1}\left( M\right) $ be the $g$-dual of $P$ i.e. $g\left( P,X\right) =\psi
\left( X\right) $ for every $X\in {\cal X}\left( M\right) $. A
straightforward computation gives:

{\bf Proposition 2.1} {\it The triple }$\left( M,g,\overline{g}\right) ${\it %
\ yields a Frobenius structure if and only if for every} $X,Y,Z\in {\cal X}%
\left( M\right) $:
$$
\theta \left( X\right) g\left( Y,Z\right) +\psi \left( Z\right) g\left(
Y,X\right) =\theta \left( Z\right) g\left( Y,X\right) +\psi \left( X\right)
g\left( Y,Z\right) .\eqno\left( 2.2\right)
$$

Let us consider the particular case when $\theta $ and $\psi $ are
proportional i.e. there exists $f\in C^{\infty }\left( M\right) $ such that $%
\psi =f\theta $. Relation $\left( 2.2\right) $ becomes:
$$
\left( 1-f\right) \theta \left( X\right) g\left( Y,Z\right) =\left(
1-f\right) \theta \left( Z\right) g\left( Y,X\right) .\eqno\left( 2.3\right)
$$

{\bf Proposition 2.2} {\it Suppose that }$f\neq 1$, $M${\it \ is connected
with }$n\geq 2${\it \ and the triple }$\left( M,g,\overline{g}\right) ${\it %
\ yields a Frobenius structure. Then }$\theta =0=\psi ,P=0$ {\it and the
Frobenius structure is degenerate i.e. }$A=0${\it .}

{\bf Proof} Let $\{U_{i}\}_{1\leq i\leq n}$ be a $g$-orthonormal basis for $%
{\cal X}\left( M\right) $. For $X=U_{i},Y=Z=U_{j},i\neq j,$ (2.3) reads $%
\theta \left( U_{i}\right) =0$. But $i$ is arbitrary. \qquad $\Box $

{\bf Particular cases}:

1) If $P=0$ then $g,\overline{g}$ are in {\it a geodesic} (or {\it projective%
}) correspondence. By a well-known result of Weyl, in this case $g$ and $%
\overline{g}$ have the same geodesics.

2) For a {\it conformal change} $\overline{g}=e^{2u}g$ with $u\in C^{\infty
}\left( M\right) $, we have $\left( 2.1\right) $ with $\theta =du$ (the
differential of $u$) and $P=-\nabla u$ (the gradient of $u$).

{\bf Application: Two-dimensional Einstein spaces}

Let $\left( M,g\right) $ be an Einstein space with $R$ the Ricci tensor;
hence $R=\lambda g$ with $\lambda \in C^{\infty }\left( M\right) \backslash
\{0\}$. For $n\geq 3$ $\lambda $ is a constant and then we confine ourselves
to the two-dimensional case. Suppose that $R$ is nondegenerate and let $%
\overline{\nabla }$ be the Levi-Civita connection of $R$. Supposing that $A$
yields a Frobenius structure we get that $\lambda $ is constant and $A=0$.
\qquad $\square $

In conclusion, returning to the general framework, the only favorable case $%
\left( 2.1\right) $ is when $\theta $ is exactly the $g$-dual of $P$. Using
a standard notation in Riemannian geometry the relation $\left( 2.1\right) $
reads:
$$
A=\theta \otimes \delta +\delta \otimes \theta +g\otimes \theta ^{\#}.\eqno%
\left( 2.1^{\prime }\right)
$$

{\bf Example 2.2 Hypersurfaces with nondegenerate second fundamental form}

Let $M$ be a hypersurface in ${I\!\!R}^{n+1}$ with $n\geq 2$ and $g,b$ the
first and second fundamental form of $M$. Suppose that $rank$ $b=n$ let $%
\nabla ,\overline{\nabla }$ be the Levi-Civita connection of $g,b$. The
corresponding $A$ it given by(\cite{s:u}):
$$
b\left( A\left( X,Y\right) ,Z\right) =-\frac{1}{2}\left( \nabla _{X}b\right)
\left( Y,Z\right) \eqno\left( 2.4\right)
$$
which yields the Frobenius relation
$$
(\nabla _{X}b)\left( Y,Z\right) =\left( \nabla _{Y}b\right) \left(
Z,X\right) .\eqno\left( 2.5\right)
$$
But this is exactly the Codazzi equation and then:

{\bf Proposition 2.3} {\it The triple }$\left( M,b,A\right) ${\it \ is a
formal Frobenius structure.}

{\bf Example 2.3 Riemannian metrics related by a self-adjoint operator}

Let $g,\widetilde{g}$ be two Riemannian metrics on $M$. Then there exists a
unique $J\in T_{1}^{1}\left( M\right) $ such that:
$$
\widetilde{g}\left( X,Y\right) =g\left( X,JY\right) =g\left( JX,Y\right) %
\eqno\left( 2.6\right)
$$
for every $X,Y,Z\in {\cal X}\left( M\right) $. If $\nabla ,\widetilde{\nabla
}$ is the Levi-Civita connection of $g,\widetilde{g}$ then:
\[
\widetilde{g}\left( A\left( X,Y\right) ,Z\right) =g\left( X,\left( \nabla
_{Y}J\right) z-\left( \nabla _{Z}J\right) X\right) +
\]
$$
+g\left( Y,\left( \nabla _{X}J\right) Z-\left( \nabla _{Z}J\right) X\right)
+g\left( Z,\left( \nabla _{X}J\right) Y+\left( \nabla _{Y}J\right) X\right) %
\eqno\left( 2.7\right)
$$
where, as usual, $A=\widetilde{\nabla }-\nabla $. Then we get:

{\bf Proposition 2.4} {\it The triple }$\left( M,\widetilde{g},A\right) $%
{\it \ is a Frobenius structure if and only if}:
$$
g\left( Y,\left( \nabla _{X}J\right) Z-\left( \nabla _{Z}J\right) X\right)
=g\left( X,\left( \nabla _{Z}J\right) Y\right) -g\left( Z,\left( \nabla
_{X}J\right) Y\right) .\eqno\left( 2.8\right)
$$

{\bf Particular case}

Let us suppose that $J$ is $\nabla ${\it -recurrent }i.e. there exists $%
\omega \in \Omega ^{1}\left( M\right) $ such that $\nabla _{X}J=\omega
\left( X\right) J$. Then:

{\bf Proposition 2.5} {\it With above condition the triple }$\left( M,%
\widetilde{g},A\right) $ {\it is a Frobenius structure if and only if}:
$$
\omega \left( X\right) g\left( Y,Z\right) =\omega \left( Z\right) g\left(
Y,X\right) .\eqno\left( 2.9\right)
$$

\section{Weak Frobenius structures generated by deformation algebras}

{\bf Example 3.1 Golab connections}

Let $\theta \in \Omega ^{1}\left( M\right) $ and $F\in T_{1}^{1}\left(
M\right) $.

{\bf Definition 3.1} A linear connection $\overline{\nabla }$ associated to
the triple $\left( g,\theta ,F\right) $ is called a {\it Golab connection}
if(\cite{g:s}):

(i) it is a metric connection i.e. $\overline{\nabla }_{X}g=0$

(ii) it has the torsion:
$$
\overline{T}\left( X,Y\right) =\theta \left( Y\right) F\left( X\right)
-\theta \left( X\right) F\left( Y\right) \text{.}\eqno\left( 3.1\right)
$$

The Golab connection exists and is unique, with expression(\cite{m:p}):
$$
\overline{\nabla }_{X}Y=\nabla _{X}Y+\theta \left( Y\right) F\left( X\right)
-S\left( X,Y\right) P\eqno\left( 3.2\right)
$$
where $S\left( X,Y\right) =g\left( FX,Y\right) $ and $P$ is the $g$-dual of $%
\theta $ i.e. $g\left( P,X\right) =\theta \left( X\right) $ for every $X\in
{\cal X}\left( M\right) $. It follows:
\[
A\left( X,Y\right) =\theta \left( Y\right) F\left( X\right) -S\left(
X,Y\right) P
\]
and a straightforward computation gives:

{\bf Proposition 3.2} {\it The Golab triple }$\left( g,\theta ,F\right) $%
{\it \ yields a weak Frobenius structure if and only if for every} $X,Y,Z\in
{\cal X}\left( M\right) :$
$$
\theta \left( X\right) S\left( Y,Z\right) +\theta \left( Y\right) S\left(
X,Z\right) =\theta \left( Z\right) \left[ S\left( X,Y\right) +S\left(
Y,X\right) \right] .\eqno\left( 3.3\right)
$$

{\bf Application: }$\lambda ${\bf -Hermitian metrics with respect to a }$%
\varepsilon ${\bf -structure}

Let $\lambda ,\varepsilon \in {I\!\!R\backslash \{0\}}$. $F$ is called $%
\varepsilon -${\it structure} if $F^{2}=\varepsilon 1_{{\cal X}\left(
M\right) }$ and $g$ is called $\lambda $-{\it Hermitian w.r.t}. $F$ if:
$$
g\left( FX,FY\right) =\lambda g\left( X,Y\right) .\eqno\left( 3.4\right)
$$

{\bf Proposition 3.3} {\it If }$g$ {\it is }$\lambda -${\it Hermitian w.r.t.
}$\varepsilon $-{\it structure }$F${\it \ and }$\left( g,\theta ,F\right) $
{\it yields a weak Frobenius structure then }$\theta =0$ {\it and therefore }%
$A=0$.

{\bf Proof }Obviously we can restrict to cases $\varepsilon ,\lambda =\pm 1$
and $\theta \left( P\right) =0$ implies $\theta =0$ because $\theta \left(
P\right) =g\left( P,P\right) $ and then $P=0$. With $Z=P$ in $\left(
3.3\right) $ we get:
$$
\theta \left( X\right) \theta \left( FY\right) +\theta \left( Y\right)
\theta \left( FX\right) =\theta \left( P\right) \left[ g\left( FX,Y\right)
+g\left( X,FY\right) \right] \eqno\left( 3.5\right)
$$
and $Y\rightarrow FY$ in $\left( 3.5\right) $ yields:
$$
\varepsilon \theta \left( X\right) \theta \left( Y\right) +\theta \left(
FX\right) \theta \left( FY\right) =\theta \left( P\right) g\left( X,Y\right)
\left( \varepsilon +\lambda \right) .\eqno\left( 3.6\right)
$$
Again $Y=P$ in $\left( 3.6\right) $ leads to:
$$
\theta \left( FX\right) \theta \left( FP\right) =\lambda \theta \left(
X\right) \theta \left( P\right) \text{.}\eqno\left( 3.7\right)
$$
With $X=P$ in $\left( 3.7\right) $:
$$
\theta \left( FP\right) ^{2}=\lambda \theta \left( P\right) ^{2}\eqno\left(
3.8\right)
$$
and then we have the cases:

I) $\lambda =-1$ implies $\theta \left( P\right) =0$.

II) $\lambda =+1.$ Relations $\left( 3.7\right) ,\left( 3.8\right) $ with $%
\lambda =1$ implies $\theta \left( FX\right) =\mu \theta \left( X\right) $
with $\mu =\pm 1$ for every $X\in {\cal X}\left( M\right) $. Plugging in $%
\left( 3.6\right) $:
$$
\left( \varepsilon +1\right) \theta \left( X\right) \theta \left( Y\right)
=\theta \left( P\right) g\left( X,Y\right) \left( \varepsilon +1\right) \eqno%
\left( 3.9\right)
$$
and we have the subclasses:

1) $\varepsilon =-1$. From $\left( 3.4\right) $ with $\lambda =1$ and $%
Y\rightarrow FY$ we get $-g\left( FX,Y\right) =g\left( X,FY\right) $ and
thus $\left( 3.5\right) $ reads $2\mu \theta \left( X\right) \theta \left(
Y\right) =0$ for every $X,Y\in {\cal X}\left( M\right) $ i.e. $\theta =0$.

2) $\varepsilon =+1$. Relation $\left( 3.9\right) $ becomes:
$$
\theta \left( X\right) \theta \left( Y\right) =\theta \left( P\right)
g\left( X,Y\right) \eqno\left( 3.10\right)
$$
for every $X,Y\in {\cal X}\left( M\right) $. Let $\{U_{i}\}_{1\leq i\leq n}$
be an $g$-orthonormal basis for ${\cal X}\left( M\right) $. The choice $%
X=U_{i},Y=U_{j},i\neq j$, in $\left( 3.10\right) $ implies $\theta \left(
U_{i}\right) \theta \left( U_{j}\right) =0$ and suppose $\theta \left(
U_{i}\right) =0$. With $X=Y=U_{i}$ in $\left( 3.10\right) $ we have $\theta
\left( U_{i}\right) ^{2}=0=\theta \left( P\right) $. \qquad $\Box $

{\bf Remarks}

(i) $\varepsilon =-1,\lambda =+1\left( -1\right) $ means that $g$ is
Hermitian(anti-Hermitian) w.r.t. {\it almost complex structure} $F$.

(ii) $\varepsilon =1,\lambda =-1$ means that $g$ is para-Hermitian w.r.t.
{\it almost product structure} $F$(\cite[p. 91]{c:fg}).

(iii) An example for $\varepsilon =\lambda =1$ is given by $F=1_{{\cal X}%
\left( M\right) }$. The Golab connection for $\left( \theta ,1_{{\cal X}%
\left( M\right) }\right) $ is called {\it Lyra} {\it connection}(\cite{l:g}).

{\bf Example 3.2 Cartan-Schouten connections on Lie groups}

Let $M=G$ be a Lie group and $\{E_{i}\}_{1\leq i\leq n}$ a basis in $L\left(
G\right) $ the Lie algebra of $G$. The {\it Cartan-Schouten connections} $%
\overline{\nabla },\stackrel{+}{\nabla },\stackrel{\circ }{\nabla }$ on $G$
are:
$$
\overline{\nabla }_{E_{i}}E_{j}=0,\quad \stackrel{+}{\nabla }_{E_{i}}E_{j}=%
\left[ E_{i},E_{j}\right] ,\quad \stackrel{\circ }{\nabla }_{E_{i}}E_{j}=%
\frac{1}{2}\left[ E_{i},E_{j}\right] \eqno\left( 3.11\right)
$$

Let $g$ be the Riemannian metric with respect to which the given basis is
orthonormal: $g\left( E_{i},E_{j}\right) =\delta _{ij}$ and let $\nabla $
the Levi-Civita connection of $g$. We define:
$$
\overline{A}=\overline{\nabla }-\stackrel{\circ }{\nabla },\stackrel{+}{A}=%
\stackrel{+}{\nabla }-\stackrel{\circ }{\nabla },A=\stackrel{+}{\nabla }-%
\overline{\nabla },\stackrel{+}{A}^{\prime }=\stackrel{+}{\nabla }-\nabla
,A^{\prime }=\stackrel{\circ }{\nabla }-\nabla ,\overline{A}^{\prime }=%
\overline{\nabla }-\nabla \eqno\left( 3.12\right)
$$
which yields(\cite[p. 35]{n:l}):
$$
\left\{
\begin{array}{l}
A\left( E_{i},E_{j}\right) =2\stackrel{+}{A}\left( E_{i},E_{j}\right) =-2%
\overline{A}\left( E_{i},E_{j}\right) =\left[ E_{i},E_{j}\right] \\
2g\left( \overline{A}^{\prime }\left( E_{i},E_{j}\right) ,E_{k}\right)
=g\left( E_{j},\left[ E_{i},E_{k}\right] \right) -g\left( E_{k},\left[
E_{i},E_{j}\right] \right) +g\left( E_{i},\left[ E_{j},E_{k}\right] \right)
\\
2g\left( \stackrel{+}{A}^{\prime }\left( E_{i},E_{j}\right) ,E_{k}\right)
=g\left( E_{i},\left[ E_{j},E_{k}\right] \right) +g\left( E_{j},\left[
E_{i},E_{k}\right] \right) +g\left( E_{k},\left[ E_{i},E_{j}\right] \right)
\\
2g\left( A^{\prime }\left( E_{i},E_{j}\right) ,E_{k}\right) =g\left( E_{i},
\left[ E_{j},E_{k}\right] \right) +g\left( E_{j},\left[ E_{k},E_{i}\right]
\right)
\end{array}
\right. \eqno\left( 3.13\right)
$$
and then:

{\bf Proposition 3.4} {\it The Cartan-Schouten triples } \newline
$\left( G,g,A\right) ,\left( G,g,\stackrel{+}{A}\right) ,\left( G,g,%
\overline{A}\right) ,\left( G,g,\overline{A}^{\prime }\right) ,\left( G,g,%
\stackrel{+}{A}^{\prime }\right) ,\left( G,g,A^{\prime }\right) ${\it \ are
weak Frobenius structures if and only if:}
$$
g\left( E_{i},\left[ E_{j},E_{k}\right] \right) =g\left( E_{j},\left[
E_{k},E_{i}\right] \right) .\eqno\left( 3.14\right)
$$

{\bf Remarks 3.5} (i) The torsions of the Cartan-Schouten connections are:
$$
\overline{T}\left( X,Y\right) =-\left[ X,Y\right] ,\quad \stackrel{+}{T}%
\left( X,Y\right) =\left[ X,Y\right] ,\quad \stackrel{\circ }{T}\left(
X,Y\right) =0\eqno\left( 3.15\right)
$$
and then $A^{\prime }$ yields exactly a Frobenius structure if $\left(
3.14\right) $ holds.

(ii) A pair $\left( {\frak g},<,>\right) $ with $\left( {\frak g},\left[ ,%
\right] \right) $ a Lie algebra and $<,>$ a scalar product on ${\frak g}$
such that:
$$
<\left[ x,y\right] ,z>+<y,\left[ x,z\right] >=0\eqno\left( 3.16\right)
$$
for every $x,y,z\in {\frak g}$ is called {\it orthogonal Lie algebra}. On a
basis $\{E_{i}\}$ of ${\frak g}$ the relation $\left( 3.16\right) $ reads
exactly as $\left( 3.14\right) $ and so we can restate the last proposition:

{\bf Proposition 3.4'} {\it The Cartan-Schouten triples are weak Frobenius
structures if and only if }$\left( {\cal X}\left( G\right) ,\left[ ,\right]
\right) ${\it \ is an orthogonal Lie algebra.}

Recall also:

{\bf Proposition 3.6} {\it The following are equivalent:}

{\it 1) the }$A${\it -algebra is commutative}

{\it 2) the }$\stackrel{+}{A}${\it -algebra is commutative}

{\it 3) the }$\overline{A}${\it -algebra is commutative}

{\it 4) the }$\overline{A}^{\prime }${\it -algebra is commutative}

{\it 5) the Lie group }$G${\it \ is abelian}

{\it 6) the Lie algebra }$L\left( G\right) ${\it \ is commutative}.

Therefore if $G$ is abelian and $\left( 3.14\right) $ holds then $A,%
\stackrel{+}{A},\overline{A},\overline{A}^{\prime }$ yields Frobenius
structures.

{\bf Example 3.3 Complex connections on K\"{a}hler manifolds}

Let $\left( M,g,J\right) $ be a K\"{a}hler manifold. Recall that a linear
connection $\overline{\nabla }$ on $M$ is called {\it complex connection} if
$\overline{\nabla }J=0$ and the family of linear connections is given by:
$$
\overline{\nabla }_{X}Y=\nabla _{X}Y+\frac{1}{2}\left( \nabla _{X}J\right)
JY+\frac{1}{2}\left( Q\left( X,Y\right) -JQ\left( X,JY\right) \right) \eqno%
\left( 3.17\right)
$$
where $\nabla $ is an arbitrary linear connection and also $Q\in
T_{2}^{1}\left( M\right) $ is arbitrary.

Set $\nabla $ the Levi-Civita connection of $g$; it results that $\nabla $
is a complex connection because $\left( M,g,J\right) $ is K\"{a}hler.
Therefore $A=\overline{\nabla }-\nabla $ is:
$$
A\left( X,Y\right) =\frac{1}{2}\left( Q\left( X,Y\right) -JQ\left(
X,JY\right) \right) \eqno\left( 3.18\right)
$$
and then:

{\bf Proposition 3.7} {\it The K\"{a}hler triple }$\left( M,g,J\right) ${\it %
\ yields a weak Frobenius structure via }$Q\in T_{2}^{1}\left( M\right) $
{\it if and only if }$Q${\it \ satisfy for every} $X,Y,Z\in {\cal X}\left(
M\right) :$%
$$
g\left( Q\left( X,Y\right) ,Z\right) -g\left( Q\left( Y,Z\right) ,X\right)
=g\left( Q\left( Y,JZ\right) ,JX\right) -g\left( Q\left( X,JY\right)
,JZ\right) .\eqno\left( 3.19\right)
$$

{\bf Example 3.4 Chern and Bismut connections on Hermitian manifolds}

Let $\left( M,g,J\right) $ be a $2n$-dimensional ($n>1$) Hermitian manifold
with complex structure $J$ and compatible metric $g$. Let $\Omega \in \Omega
^{2}\left( M\right) $ be {\it the K\"{a}hler form}, $\Omega \left( \cdot
,\cdot \right) =g\left( \cdot ,J\cdot \right) $, and $\theta \in \Omega
^{1}\left( M\right) $ {\it the Lee form}, $\theta =\frac{1}{n-1}d^{*}\Omega
\circ J$. Let $\nabla $ be the Levi-Civita connection of $g$ and $\nabla
^{C},\nabla ^{B}$ the {\it Chern}, respectively {\it Bismut}, {\it connection%
} on $\left( M,g,J\right) $:
$$
\left\{
\begin{array}{l}
g\left( \nabla _{X}^{C}Y,Z\right) =g\left( \nabla _{X}Y,Z\right) +\frac{1}{2}%
d\Omega \left( JX,Y,Z\right) \\
g\left( \nabla _{X}^{B}Y,Z\right) =g\left( \nabla _{X}Y,Z\right) -\frac{1}{2}%
d\Omega \left( JX,JY,JZ\right)
\end{array}
\right. .\eqno\left( 3.20\right)
$$
For remarkable properties of these connections see \cite{g:h}.

Therefore, denoting $A^{C}=\nabla ^{C}-\nabla ,A^{B}=\nabla ^{B}-\nabla $
the Chern, respectively Bismut, deformation tensor it results:
$$
\left\{
\begin{array}{l}
g\left( A^{C}\left( X,Y\right) ,Z\right) =\frac{1}{2}d\Omega \left(
JX,Y,Z\right) \\
g\left( A^{B}\left( X,Y\right) ,Z\right) =-\frac{1}{2}d\Omega \left(
JX,JY,JZ\right)
\end{array}
\right. \eqno\left( 3.21\right)
$$
and a straightforward computation gives:

{\bf Proposition 3.8} (i) {\it The triple }$\left( M,g,A^{C}\right) ${\it \
is a weak Frobenius structure if and only if}:
$$
d\Omega \left( JX,Y,Z\right) =d\Omega \left( JY,Z,X\right) .\eqno\left(
3.22\right)
$$
(ii) {\it The triple }$\left( M,g,A^{B}\right) ${\it \ is a weak Frobenius
structure if and only if the 3-form }$d\Omega \left( J\cdot ,J\cdot ,J\cdot
\right) ${\it \ is invariant to circular permutations i.e}.
$$
d\Omega \left( JX,JY,JZ\right) =d\Omega \left( JY,JZ,JX\right) .\eqno\left(
3.23\right)
$$

{\bf Particular case} $d\Omega =\theta \wedge \Omega $

For $n=2$ i.e. in the 4-dimensional case or $n>2$ and $\left( M,g,J\right) $
is {\it locally conformal K\"{a}hler} (another notion introduced by I.
Vaisman, \cite{v:l}) we have $d\Omega =\theta \wedge \Omega $, \cite{d:o}.
Then we obtain:

{\bf Proposition 3.9} {\it Let }$\left( M,g,J\right) ${\it \ as above.}

{\it (i) The triple }$\left( M,g,A^{C}\right) ${\it \ is a weak Frobenius
structure if and only if}:
\[
d^{*}\Omega \left( X\right) \Omega \left( Y,Z\right) +d^{*}\Omega \left(
JY\right) g\left( Z,X\right) -d^{*}\Omega \left( JZ\right) g\left(
X,Y\right) =
\]
$$
=d^{*}\Omega \left( Y\right) \Omega \left( Z,X\right) +d^{*}\Omega \left(
JZ\right) g\left( X,Y\right) -d^{*}\Omega \left( JX\right) g\left(
Y,Z\right) .\eqno\left( 3.24\right)
$$
(ii) {\it The triple }$\left( M,g,A^{B}\right) ${\it \ is a weak Frobenius
structure if and only if the 3-form }$d^{*}\Omega \left( \cdot \right)
\Omega \left( \cdot ,\cdot \right) ${\it \ is invariant to circular
permutations}:
\[
d^{*}\Omega \left( X\right) \Omega \left( Y,Z\right) +d^{*}\Omega \left(
Y\right) \Omega \left( Z,X\right) +d^{*}\Omega \left( Z\right) \Omega \left(
X,Y\right) =
\]
$$
=d^{*}\Omega \left( Y\right) \Omega \left( Z,X\right) +d^{*}\Omega \left(
Z\right) \Omega \left( X,Y\right) +d^{*}\Omega \left( X\right) \Omega \left(
Y,Z\right) .\eqno\left( 3.25\right)
$$

\section{(Weak) Frobenius structures with $A$ not of $\left( \overline{%
\protect\nabla }-\protect\nabla \right) $-type}

''{\bf Pseudo-''Example 4.1} {\bf Cross products on }${I\!\!R}^{3}$ and ${%
I\!\!R}^{7}$

An immediately example of skew-symmetric weak Frobenius structure is given
by: $M={I\!\!R}^{3}$ or ${I\!\!R}^{7},g=<,>$ the Euclidean inner product and
$A=\times $ the usual cross product, because for every $x,y,z\in {I\!\!R}%
^{3} $ we have $<x,y\times z>=<y,z\times x>=<z,x\times y>$. In ${I\!\!R}^{3}$
this equality is well-known and for ${I\!\!R}^{7}$ see, for example,
\cite[p. 71]{p:y}. But the cross product is a pseudo-tensor and not a tensor!

{\bf Example 4.2 Frobenius structures generated by selfadjoint o\-pe\-rators}

Fix $J\in T_{1}^{1}\left( M\right) $ which is $g$-{\it selfadjoint} i.e. for
every $X,Y\in X\left( M\right) $:
$$
g\left( JX,Y\right) =g\left( X,JY\right) .\eqno\left( 4.1\right)
$$

{\bf Proposition 4.1} {\it If }$J${\it \ is self-adjoint and for every }$%
X,Y\in X\left( M\right) ${\it :}
$$
\left( \nabla _{X}J\right) Y=\left( \nabla _{Y}J\right) X\eqno\left(
4.2\right)
$$
{\it then }$\left( M,g,A\right) ${\it \ with }$A\left( X,Y\right) =\left(
\nabla _{X}J\right) Y${\it \ is a Frobenius structure}.

{\bf Proof} Let us apply $\nabla _{X}$ to $\left( 4.1\right) $ written in $Y$
and $Z$:
\[
\nabla _{X}\left( g\left( JY,Z\right) \right) =\nabla _{X}\left( g\left(
Y,JZ\right) \right)
\]
which yields:
\[
g\left( \nabla _{X}JY,Z\right) +g\left( Y,J\left( \nabla _{X}Z\right)
\right) =g\left( \nabla _{X}JZ,Y\right) +g\left( J\left( \nabla _{X}Y\right)
,Z\right)
\]
or:
\[
g\left( \nabla _{X}JY-J\left( \nabla _{X}Y\right) ,Z\right) =g\left(
Y,\nabla _{X}JZ-J\left( \nabla _{X}Z\right) \right)
\]
which means:
\[
g\left( \left( \nabla _{X}J\right) Y,Z\right) =g\left( Y,\left( \nabla
_{X}J\right) Z\right) \stackrel{\left( 4.2\right) }{=}g\left( Y,\left(
\nabla _{Z}J\right) X\right)
\]
a relation equivalent with Frobenius condition $\left( 1.1\right) $. Also, $%
\left( 4.2\right) $ means the commutativity of $A$. \qquad $\Box $

{\bf Particular cases}:

1) Let $\left( M,g\right) $ be an {\it invariant} hypersurface of a
Riemannian manifold i.e. the curvature tensor of $M$ is tangent to $M$. Let $%
J$ be the Weingarten operator of $M$. Then $\left( 4.2\right) $ is exactly
the Codazzi equation and $\left( 4.1\right) $ holds because:
\[
g\left( JX,Y\right) =b\left( X,Y\right) =b\left( Y,X\right) =g\left(
X,JY\right)
\]
where $b$ is the second fundamental form of $M$.

2) Let $R\in T_{2}^{0}\left( M\right) $ be the Ricci tensor of $\left(
0,2\right) $-type and $ric\in T_{1}^{1}\left( M\right) $ defined by:
$$
g\left( ricX,Y\right) =R\left( X,Y\right) \eqno\left( 4.3\right)
$$
for every $X,Y\in {\cal X}\left( M\right) $. Then $J=ric$ is self-adjoint
and:

{\bf Proposition 4.2} {\it If for every }$X,Y,Z\in {\cal X}\left( M\right) $%
{\it :}
$$
\left( \nabla _{X}R\right) \left( Y,Z\right) =\left( \nabla _{Y}R\right)
\left( X,Z\right) \eqno\left( 4.4\right)
$$
{\it then }$\left( 4.2\right) ${\it \ holds for }$J=ric$.

{\bf Proof} From $\left( 4.4\right) $:
\[
\left( \nabla _{X}R\right) \left( Y,Z\right) -\left( \nabla _{Y}R\right)
\left( X,Z\right) =0=
\]
\[
=X\left( R\left( Y,Z\right) \right) -Y\left( R\left( X,Z\right) \right)
+R\left( X,\nabla _{Y}Z\right) -R\left( \nabla _{X}Z,Y\right) -R\left( \left[
X,Y\right] ,Z\right)
\]
and then:
$$
R\left( \left[ X,Y\right] ,Z\right) =X\left( R\left( Y,Z\right) \right)
-Y\left( R\left( X,Z\right) \right) +R\left( X,\nabla _{Y}Z\right) -R\left(
\nabla _{X}Z,Y\right) .\eqno\left( 4.5\right)
$$
But:
\[
\left( \nabla _{X}ric\right) Y-\left( \nabla _{Y}ric\right) X=\nabla
_{X}ricY-\nabla _{Y}ricX-ric\left( \left[ X,Y\right] \right)
\]
and therefore:
\[
g\left( \left( \nabla _{X}ric\right) Y-\left( \nabla _{Y}ric\right)
X,Z\right) =
\]
\[
=g\left( \nabla _{X}ricY,Z\right) -g\left( \nabla _{Y}ricX,Z\right) -g\left(
ric\left[ X,Y\right] ,Z\right) =
\]
\[
=X\left( g\left( ricY,Z\right) \right) -Y\left( g\left( ricX,Z\right)
\right) -g\left( ricY,\nabla _{X}Z\right) +
\]
\[
+g\left( ricX,\nabla _{Y}Z\right) +R\left( \left[ X,Y\right] ,Z\right)
\stackrel{\left( 4.5\right) }{=}0.
\]
Because $Z\in {\cal X}\left( M\right) $ is arbitrary we get the conclusion.
\qquad $\Box $

\section{Examples of symmetric non-metric linear connections}

Necessary and sufficient conditions for a symmetric linear connection to be
Levi-Civita with respect to a Riemannian metric are given in \cite{s:c} and
two examples of non-metric connections are presented in the cited paper. In
this section we give other two examples using the differential system of
autoparallel curves. It is amazing that all three papers cited in this
section are from physics oriented journals (Comm. Math. Phys., Rep. Math.
Phys., J. Phys. A) not from pure mathematical journals!

{\bf Example 5.1 }(two-dimensional)

In \cite[p. 100-101]{b:s} it is proved that the differential system:
$$
\left\{
\begin{array}{l}
\stackrel{..}{x}+\frac{y}{1+y^{2}}\stackrel{.}{x}\stackrel{.}{y}=0 \\
\stackrel{..}{y}=0
\end{array}
\right. \eqno\left( 5.1\right)
$$
may be naturally interpreted as the system of autoparallel curves for a
torsion-free connection which is not a metric connection( two proofs are
given). The first integrals of $\left( 5.1\right) $ are (see the arguments
of next example):
$$
F_{1}=\stackrel{.}{y},\quad F_{2}=\stackrel{.}{x}\sqrt{1+y^{2}}.\eqno\left(
5.2\right)
$$

{\bf Example 5.2} (three-dimensional)

In \cite[p. 2188]{t:p} it is proved that {\it the Halphen system}:
$$
\left\{
\begin{array}{l}
\stackrel{.}{x}_{1}=x_{2}x_{3}-x_{3}x_{1}-x_{1}x_{2} \\
\stackrel{.}{x}_{2}=x_{3}x_{1}-x_{1}x_{2}-x_{2}x_{3} \\
\stackrel{.}{x}_{3}=x_{1}x_{2}-x_{2}x_{3}-x_{3}x_{1}
\end{array}
\right. \eqno\left( 5.3\right)
$$
admits no polynomial first integrals. Therefore, with the change $%
x_{i}\rightarrow \stackrel{.}{x}_{i}$ we conclude that the differential
system:
$$
\left\{
\begin{array}{l}
\stackrel{..}{x}_{1}=\stackrel{.}{x}_{2}\stackrel{.}{x}_{3}-\stackrel{.}{x}%
_{3}\stackrel{.}{x}_{1}-\stackrel{.}{x}_{1}\stackrel{.}{x}_{2} \\
\stackrel{..}{x}_{2}=\stackrel{.}{x}_{3}\stackrel{.}{x}_{1}-\stackrel{.}{x}%
_{1}\stackrel{.}{x}_{2}-\stackrel{.}{x}_{2}\stackrel{.}{x}_{3} \\
\stackrel{..}{x}_{3}=\stackrel{.}{x}_{1}\stackrel{.}{x}_{2}-\stackrel{.}{x}%
_{2}\stackrel{.}{x}_{3}-\stackrel{.}{x}_{3}\stackrel{.}{x}_{1}
\end{array}
\right. \eqno\left( 5.4\right)
$$
does not admit a polynomial first integral. If this symmetric connection
will be Levi-Civita for the metric $g=\left( g_{ij}\right) $ then the
Hamiltonian $H=\frac{1}{2}g_{ij}\stackrel{.}{x}_{i}\stackrel{.}{x}_{j}$,
which is 2-homogeneous, will be a first integral for $\left( 5.4\right) $,
false.

\vspace{.3cm}

\noindent Faculty of Mathematics \newline
University "Al. I. Cuza" \newline
Ia\c si, 6600 \newline
Romania \newline
email: mcrasm@uaic.ro

\end{document}